\magnification=1200
\input amstex
\documentstyle{amsppt}

\def\wt#1{\widetilde {#1}}
\def\ssm {\smallsetminus}
\def\SC {\Cal C}
\def\SD {\Cal D}
\def\SE {\Cal E}
\def\SF {\Cal F}
\def\SG {\Cal G}
\def\SH {\Cal H}
\def\SK {\Cal K}
\def\SL {\Cal L}
\def\SN {\Cal N}
\def\SO {\Cal O}
\def\SP {\Cal P}
\def\SQ {\Cal Q}
\def\SR {\Cal R}
\def\SU {\Cal U}
\def\SW {\Cal W}
\def\SX {\Cal X}
\def\SY {\Cal Y}

\topmatter
\title Moduli spaces of $\text{SL}(r)$-bundles on singular 
irreducible curves
 \endtitle
\author Xiaotao Sun   \endauthor
\address Institute of Mathematics, Chinese Academy of Sciences,
Beijing 100080, China \endaddress \email
xsun$\@$math08.math.ac.cn\endemail
\address Department of Mathematics, The University of
Hong Kong, Pokfulam Road, Hong Kong\endaddress \email
xsun$\@$maths.hku.hk\endemail
\thanks The work is supported by a grant of NFSC for outstanding
young researcher at contract number 10025103 and a grant of
RGC of Hong Kong at RGC code 7131/02P.
\endthanks
\endtopmatter
\document

\heading{Introduction}\endheading

One of the problems in moduli theory, motivated by physics, is
to study the degeneration of moduli spaces of semistable
$G$-bundles on curves of genus $g\ge 2$.
When a smooth curve $Y$ specializes to a stable curve $X$,
one expects that the moduli space of semistable $G$-bundles on $Y$
specializes
to a (nice) moduli space of generalized semistable $G$-torsors on $X$.
It is well known ([Si]) that for
any flat family $\SC\to S$ of stable curves 
there is a family $\SU(r,d)_S\to S$ of moduli
spaces $\SU_{\SC_s}(r,d)$ of
(s-equivalence classes of) semistable torsion free sheaves of rank $r$
and degree $d$ on curves $\SC_s$ ($s\in S$).
If we fix a suitable representation $G\to \text{GL}(r)$, one would like
to define a closed subscheme $\SU_X(G)\subset\SU_X(r,d)$ 
which should be a moduli space of sutiable $G$-sheaves on $X$. Moreover,
it should behave well under specialization, i.e. if a smooth curve $Y$
specializes to $X$, then the moduli space of $G$-bundles on $Y$ specializes
to $\SU_X(G)$.
By my knowledge, the problem is almost completely open except for special
case like $G=\text{SO}(r)$ or $G=\text{Sp}(r)$ ([Fa1], [Fa2]), 
where one has a generalisation of $G$-torsors which extends the case 
$G=\text{GL}(r)$. It is
open even for $G=\text{SL}(r)$ (See [Fa1], [Fa2] for the introduction).

In this paper, we will consider the case $G=\text{SL}(r)$ and $X$ being
irreducible (the case of a reducible curve with one node was studied in
[Su2]).  For any projective curve $X$,  we will
use $\SU_X(r,d)$ to denote the moduli space of semistable 
torsion free sheaves of rank $r$ and degree $d$ on $X$. If 
$X_{\eta}$ is a smooth curve and $L_{\eta}$ is a line bundle of degree $d$ 
on $X_{\eta}$,
we use $\SU_{X_{\eta}}(r,L_{\eta})$ to denote the moduli space 
of semistable vector 
bundles of
rank $r$ with fixed determinant $L_{\eta}$ on $X_{\eta}$, 
which is a closed subvariety
of $\SU_{X_{\eta}}(r,d)$. It is known that when $X_{\eta}$ specializes to 
$X$ the moduli
space $\SU_{X_{\eta}}(r,d)$ specializes to $\SU_X(r,d)$. It is natural to expect
that if $L_{\eta}$ specializes to a torsion free sheaf $L$ on $X$ then 
$\SU_{X_{\eta}}(r,L_{\eta})$ specializes to a closed subscheme 
$\SU_X(r,L)\subset \SU_X(r,d)$. It is important that we should look for
an intrinsic $\SU_X(r,L)$ (i.e. independent of $X_{\eta}$) which should
not be {\it too bad} and should represent a {\it moduli problem}.

Let $S=\text{Spec($A$)}$ where $A$ is a discrete valuation ring,
let $\SC\to S$ be a proper flat family of curves with closed fibre
$\SC_0\cong X$ and smooth generic fibre $\SC_{\eta}$. Then we have
a $S$-flat scheme $\SU(r,d)_S\to S$ with generic fibre 
$\SU_{\SC_{\eta}}(r,d)$ and closed fibre being $\SU_X(r,d)$.
For any line bundle
$\SL_{\eta}$ of degree $d$ on $\SC_{\eta}$, there is a unique
extension $\SL$ on $\SC$ such that $\SL|_{\SC_0}:=L$ is torsion free 
of degree $d$ (since $X$ is irreducible). Then 
$\SU_{\SC_{\eta}}(r,\SL_{\eta})\subset \SU(r,d)_S$ is an 
irreducible, reduced, locally closed subscheme. Let
$$f:\SU(r,\SL)_S:=\overline{\SU_{\SC_{\eta}}(r,\SL_{\eta})}
\subset \SU(r,d)_S\to S$$
be the Zariski closure of $\SU_{\SC_{\eta}}(r,\SL_{\eta})$ in $\SU(r,d)_S$.
Then $f:\SU(r,\SL)_S\to S$ is flat and projective, but there is
no reason that its closed fibre $f^{-1}(0)$ (even its support 
$f^{-1}(0)_{\text{red}}$) is indpendent of the family $\SC\to S$ and
$\SL_{\eta}$. However, 
there are conjectures ([NS]) that $f^{-1}(0)$ is intrinsic 
for irreducible curves $X$ with only one node. 
To state them, we introduce the notation
for any stable irreducible
curves. Let $X$ be an irreducible stable curve with $\delta$ nodes
$\{x_1,\,...,\,x_{\delta}\}$, and $L$ a torsion free sheaf of rank
one and degree $d$ on
$X$. A torsion free sheaf $F$ of rank $r$ and degree $d$
on $X$ is called with a determinant $L$ if there exists a 
morphism $(\wedge^rF)\to L$ which is an isomorphism outside the nodes
of $X$. The subset $\SU_X(r,L)\subset\SU_X(r,d)$
consists of $s$-equivalence classes $[F]\in\SU_X(r,d)$ such that
$[F]$ contains a sheaf with a fixed determinant $L$. 
Then D.S. Nagaraj and C.S. Seshadri made the following conjectures 
(See Conjecture (a) and (b) at page 136 of [NS]):
\roster 
\item If $L$ is a line bundle on $X$ and $\SU_X(r,L)^0\subset\SU_X(r,d)$
 is the subset of locally free sheaves, then
$\SU_X(r,L)$ is the closure of $\SU_X(r,L)^0$ in $\SU_X(r,d)$. 
\item
Let $\SL_{\eta}$ (resp. $L$) be a line bundle 
(resp. torsion free sheaf of rank one) of degree $d$ on smooth curve $Y$ 
(resp. $X$). 
Assume that $\SL_{\eta}$ specializes to $L$ as $Y$ specializes to $X$. 
Then $\SU_X(r,L)$ is the specialization of $\SU_Y(r,\SL_{\eta})$. 
\endroster

We answer $(1)$ completely. In fact, even if $L$ is not locally free
(thus $\SU_X(r,L)$ contains no locally free sheaf), we prove that torsion 
free sheaves of type $1$ (See Section 1) are dense in $\SU_X(r,L)$.
 
\proclaim{Theorem 1} Let $L$ be a torsion free sheaf of rank $1$ and degree
$d$. Define 
$$\SU_X(r,L)^0=\{F\in\SU_X(r,L)\,|\,(\wedge^rF)\cong L\}$$ 
which coincides with the subset of locally free sheaves when $L$
is locally free. Then 
\roster
\item $\SU_X(r,L)$ is the closure of $\SU_X(r,L)^0$. 
If $L$ is not locally free, $\SU_X(r,L)^0$ is the subset of torsion free 
sheaves of type $1$.
\item There is a canonical scheme structure on $\SU_X(r,L)^0$, which
is reduced when $L$ is locally free,
such that when smooth
curve $\SC_{\eta}$ specializes to $X$ and $\SL_{\eta}$ specializes to 
$L$ on $X$,
the specialization $f^{-1}(0)$ of $\SU_{\SC_{\eta}}(r,\SL_{\eta})$ 
contains a dense open subscheme which is isomorphic to $\SU_X(r,L)^0$.
In particular, $$f^{-1}(0)_{\text{red}}\cong\SU_X(r,L).$$
\endroster
\endproclaim

If the specialization $f^{-1}(0)$ has no $embedded$ $point$, then our
theorem also proved Conjecture (2). Unfortunately, $\SU_X(r,L)$ seems
not represent a nice moduli functor, we can not say anything about
the scheme structure of $\SU_X(r,L)$. To remedy this, 
we consider the
specialization of $\SU_{\SC_{\eta}}(r,\SL_{\eta})$ in the so called
generalized Gieseker space $G(r,d)$ (See [NSe]). Let $X$ be an irreducible
stable curve with only one node $p_0$ and $L$ be a line bundle of degrre
$d$ on $X$. 
Then, when $r\le 3$, or $r=4$ and
the normalization $\wt X$ is not hyperelliptic, we show that there is
a Cohen-Macaulay closed subscheme $G(r,L)\subset G(r,d)$ of
pure dimension $(r^2-1)(g-1)$, which represents a nice moduli functor 
(See Definition 3.2). Moreover, $G(r,L)$ satisfies the requirements in 
(2) for specializations. It is known ([NSe] that
there is a canonical birational morphism
$\theta: G(r,d)\to \SU_X(r,d)$. We prove in Lemma 3.4 that the 
{\it set-theoretic image} of $G(r,L)$ is $\SU_X(r,L)$. Thus we can endow 
$\SU_X(r,L)$ a scheme structure by the {\it scheme-theoretic image}
of $G(r,L)$. Then we have

\proclaim{Theorem 2} Let $X$ be an irreducible curve of genus $g\ge 2$
with only one node $p_0$. Let $L$ be a line bundle of degree $d$ on $X$.
Assume that $r\le 3$, or $r=4$ and the normalization of $X$ is not 
hyperelliptic. Then, when $(r,d)=1$, we have
\roster\item There is a Cohen-Macaulay projective scheme $G(r,L)$ of
pure dimension $(r^2-1)(g-1)$, which represents a moduli functor. 
\item Let $\SC\to S$ be a proper family of curves over a discrete valuation
ring, which has smooth generic fibre $\SC_{\eta}$ 
and closed fibre $\SC_0\cong X$. If there is a line bundle $\SL$ on $\SC$
such that $\SL|_{\SC_0}\cong L$. Then there exists an irreducible, reduced,
Cohen-Macaulay $S$-projective scheme 
$f: G(r,\SL)_S\to S$, which represents a
moduli functor, such that
$f^{-1}(0)\cong G(r,L),\quad 
f^{-1}(\eta)\cong\SU_{\SC_{\eta}}(r,\SL_{\eta}).$
\item There exists a proper birational $S$-morphism 
$\theta: G(r,\SL)_S\to \SU(r,\SL)_S$ which induces a birational morphism
$\theta: G(r,L)\to \SU_X(r,L).$  
\endroster
\endproclaim

 Theorem 1 is proved in Section 1. In Section 2, we introduce the objects
which are used to define Gieseker moduli space. Then Theorem 2 is proved
in Section 3.

{\it Acknowledgements.} I would like to thank Prof. C. S. Seshadri 
very much. 
Discussions and email exchanges with him are
very helpful, which stimulated
the use of generalized parabolic bundles
in Section 1 and Lemma 3.5 in Section 3.

\heading\S1 Torsion-free sheaves with fixed determinant on irreducible
curves \endheading

Let $X$ be a stable irreducible curve of genus $g$ with $\delta$ nodes
$x_1$,..., $x_{\delta}$. Any torsion free sheaf $\SF$ of rank $r$ on $X$
can be written into (locally at $x_i$)
$$\SF\otimes\hat\SO_{X,x_i}\cong
\hat\SO_{X,x_i}^{\oplus a_i}\oplus m_{x_i}^{\oplus (r-a_i)}.$$
We call that $\SF$ has type $r-a_i$ at $x_i$. 
Let $\SU_X(r,d)$ be the moduli space of
$s$-equivalence
classes of semistable torsion free sheaves of rank $r$ and degree $d$ on
$X$. Inspired by [NS], we make the following definition.

\proclaim{Definition 1.1} Let $L$ be a torsion free sheaf of rank
one and degree $d$ on $X$. A torsion free sheaf $\SF$ of rank $r$
and degree $d$ on $X$ is called with a determinant $L$ if there
exists a non-trivial morphism $\wedge^r\SF\to L$ which is an
isomorphism outside the nodes.\endproclaim

\proclaim{Lemma 1.2} For any exact sequence
$0\to\SF_1 @>\alpha>>\SF@>\beta>>\SF_2\to 0$ of torsion free sheaves
with rank $r_1$, $r$,
$r_2$ respectively, we have a morphism
$$(\wedge^{r_1}\SF_1)\otimes(\wedge^{r_2}\SF_2)\to
\frac{\wedge^r\SF}{torsion},$$
which is isomorphic outside the nodes.
In particular, if a semistable sheaf $\SF$ has a fixed determinant $L$,
then the associated graded torsion free sheaf
$gr(\SF)$ will also have the fixed determinant $L$.
\endproclaim

\demo{Proof} There is a morphism
$\wedge^{r_2}\SF_2\to \Cal Hom(\wedge^{r_1}\SF_1,\wedge^r\SF/torsion),$
which locally is defined as follows: For any $\omega\in\wedge^{r_2}\SF_2$,
choose a preimage $\wt\omega\in\wedge^{r_2}\SF$ with respect to
$\wedge^{r_2}\beta$. Then the image of $\omega$ is defined to be
the morphism $$\wedge^{r_1}\SF_1\to \wedge^r\SF/torsion,$$
which takes any $f\in\wedge^{r_1}\SF_1$ to
the section $(\wedge^{r_1}\alpha)(f)\wedge\wt\omega\in\wedge^r\SF/torsion$,
which does not depend on the choice of $\wt\omega$ since
the image of $\wedge^{r_1+1}\alpha$ is a torsion sheaf. The morphism defined above
is isomorphism outside the nodes (See Lemma 1.2 of [KW]). Thus we have
the desired morphism
$$(\wedge^{r_1}\SF_1)\otimes(\wedge^{r_2}\SF_2)\to
(\wedge^{r_1}\SF_1)\otimes\Cal Hom(\wedge^{r_1}\SF_1,
\frac{\wedge^r\SF}{torsion})\to
\frac{\wedge^r\SF}{torsion}.$$
\enddemo

\proclaim{Definition 1.3} The subset $\SU_X(r,L)\subset\SU_X(r,d)$
and $\SU_X(r,L)^0\subset\SU_X(r,L)$ are defined to be
$$\SU_X(r,L)=\left\{\aligned
&\text{$s$-equivalence classes $[\SF]\in\SU_X(r,d)$ such that}\\
&\text{$[\SF]$ contains a sheaf with a fixed determinant $L$}
\endaligned\right\}$$
$$\SU_X(r,L)^0=\left\{\,[\SF]\in\SU_X(r,L) \,|\,\wedge^r\SF\cong L\right\}$$
When $L$ is a line bundle,  $\SU_X(r,L)^0$
consists of locally free sheaves with the fixed determinant $L$.
When $L$ is not a line bundle, $\SU_X(r,L)^0$ consists of torsion
free sheaves of type $1$ at each node of $X$.\endproclaim

We first consider the case that $L$ is a line bundle and $X$ has only 
one node $p_0$. Let $\pi: \wt X\to X$ be the
normalization with $\pi^{-1}(p_0)=\{p_1,p_2\}$. The normalization
$\phi:\SP\to \SU_X(r,d)$ was studied in [Su1], where $\SP$ is the
moduli spaces of semistable generalized parabolic bundles (GPB) of
degree $d$ and rank $r$ on $\wt X$. A GPB of degree $d$ and rank
$r$ on $\wt X$ is a pair $(E,Q)$ consisting of a vector bundle $E$
of degree $d$ and rank $r$ on $\wt X$ and a $r$-dimensional
quotient $E_{p_1}\oplus E_{p_2}\to Q$. There is a flat morphism
(See Lemma 5.7 of [Su1])
$$Det:\SP\to J^d_{\wt X}$$
sending $(E,Q)$ to $det(E)$. Let $\wt L=\pi^*(L)$ and $\SP^{\wt
L}=Det^{-1}(\wt L)$. Then $\SP^{\wt L}$ is an irreducible
projective variety (See the proof of Lemma 5.7 in [Su1]).
Let $\SD_i$ ($i=1,2$) be the divisor consisting
of $(E,Q)$ such that $E_{p_i}\to Q$ is not an isomorphism (See
[Su1] for details). Let $\SD^{\wt L}_i=\SD_i\cap\SP^{\wt L}$.

\proclaim{Lemma 1.4} The set $\SU_X(r,L)$ is contained in the
image $\phi(\SP^{\wt L})$. Moreover, 
$$\SU_X(r,L)\ssm \SU_X(r,L)^0\subset \phi(\SD_1^{\wt
L}\cap\SD_2^{\wt L}).$$
\endproclaim

\demo{Proof} Let $F\in \SU_X(r,L)$ with
$F\otimes\hat\SO_{p_0}\cong{\hat\SO}^{\oplus a}_{p_0}\oplus
m_{p_0}^{\oplus(r-a)}$. Let $\wt E=\pi^*F/torsion$. Then,
by local computions (See, for example, Remark 2.1, 2.6 of [NS]), we have
$$0\to F@>d>>\pi_*\wt E\to\,_{p_0}\wt Q\to 0\tag1.1$$
where $dim(\wt Q)=a$ and the quotient $\pi_*\wt E\to\,_{p_0}\wt Q$
induces two surjective maps $\wt E_{p_i}\to \wt Q$ ($i=1,2$).
Denote their kernel by $K_i$, we have
$$0\to K_i\to\wt E_{p_i}\to\wt Q\to 0.$$
On the other hand, for $F\in \SU_X(r,L)$, 
let $\SQ$ be the cokernel of $\wedge^rF\to L$, then
$$0\to det(\wt E)\to\wt L\to\pi^*\SQ\to 0$$
where $\pi^*\SQ=\,_{p_1}V_1\oplus\,_{p_2}V_2$ and $n_1$, $n_2$ 
is respectively the dimension of $V_1$, $V_2$. 
Thus $det(\wt E)=\wt
L\otimes\SO_{\wt X}(-n_1p_1-n_2p_2)$ where $n_i\ge 0$ and
$n_1+n_2=r-a$. 

Let $h:\wt E\to E$ be the Hecke modifications at
$p_1$ and $p_2$ such that $ker(h_{p_i})\subset K_i$ has dimension
$n_i$ for $i=1,2$. Then we have
$$0\to \wt E@>h>> E\to\,_{p_1}\wt Q_1\oplus\,_{p_2}\wt Q_2\to 0\tag1.2$$
with $dim(\wt Q_i)=n_i$. Thus $det(E)=det(\wt E)\otimes\SO_{\wt
X}(n_1p_1+n_2p_2)=\wt L$ and $\phi(E,Q)=F$ if we define $Q$ by the
exact sequence
$$0\to F@>(\pi_*h)\cdot d>>\pi_*E\to\,_{p_0} Q\to 0.\tag1.3$$

To describe 
the GPB $(E,E_{p_1}\oplus E_{p_2}@>q>>Q\to 0)$, note that (1.3) induces
$$F_{p_0}@>d_{p_0}>>\wt E_{p_1}\oplus\wt E_{p_2}@>h_{p_1}\oplus
h_{p_2}>>
E_{p_1}\oplus E_{p_2}@>q>> Q\to 0.$$
Then $d_{p_0}(F_{p_0})\cap \wt E_{p_i}=K_i$ by
(1.1) and $h_{p_i}(K_i)=ker(q_i)$ by the exactness of (1.3), where
$q_i:E_{p_i}\to Q$ ($i=1,2$) are projections induced by 
$E_{p_1}\oplus E_{p_2}@>q>>Q\to 0$. Thus $dim(ker(q_i))=r-a-n_i$ by the
construction of $h$.

For any $F\in \SU_X(r,L)\ssm \SU_X(r,L)^0$, the cokernel $\SQ$ of 
$\wedge^rF\to L$ must be
non-trivial. This implies that both $V_1$ and $V_2$ in
$\pi^*\SQ=\,_{p_1}V_1\oplus\,_{p_2}V_2$
are non-trivial since
for any $i=1,2$, we have
$$Hom_{\SO_{\wt X}}(\,_{p_i}V_i,\,_{p_i}\Bbb C)=
Hom_{\SO_{\wt X}}(\pi^*\SQ,\,_{p_i}\Bbb C)=
Hom_{\SO_X}(\SQ,\pi_*(\,_{p_i}\Bbb C))\neq 0.$$
Thus their dimensions $n_1$ and $n_2$ must be positive and 
$n_1+n_2=r-a$, which means that
$ker(q_i)\neq 0$ ($i=1,2$) 
and the GPB $(E,Q)$ must be in
$\SD_1\cap\SD_2$. Thus
$$\SU_X(r,L)\ssm \SU_X(r,L)^0\subset \phi(\SD_1^{\wt
L}\cap\SD_2^{\wt L}).$$
\enddemo

\remark{Remark 1.5} 
This is also indicated in the following
consideration. There is a $\Bbb P^1$-bundle $p:\Bbb P\to J^d_{\wt
X}$ and the normalization map $\phi_1:\Bbb P\to J^d_X$. The
morphism $Det: \SP\to J^d_{\wt X}$ can be lift to a rational
morphism
$$\wt{Det}:\SP\dashrightarrow \Bbb P@>\phi_1>> J^d_X,$$
which is well-defined on $\SP\ssm \SD_1\cap\SD_2$. When $L$ is a
line bundle, $\wt{Det}^{-1}(L)$ is disjoint with
$\SD_i\ssm(\SD_1\cap\SD_2)$.\endremark

\proclaim{Lemma 1.6} Let $\Lambda$ be a discrete valuation ring
and $T=\text{Spec($\Lambda$)}$. Then, for any $F\in\SU_X(r,L)$,
there is a $T$-flat sheaf $\SF$ on $X\times T$ such that\roster
\item $\SF_t=\SF|_{X\times\{t\}}$ is locally free for $t\neq 0$ and
$\SF_0=F$,
\item  $\wedge^r(\SF|_{X\times(T\ssm\{0\})})=p_X^*L$.\endroster
In particular, $\SU_X(r,L)^0$ is dense in $\SU_X(r,L)$.
\endproclaim

\demo{Proof} Let $(E,Q)\in \SP^{\wt L}$ be the GPB such that
$\phi(E,Q)=F$ (Lemma 1.4). Then there exists a $T$-flat family of
vector bundles $\SE$ on $\wt X\times T$ with $det(\SE)=p_{\wt
X}^*\wt L$, and a $T$-flat quotient
$$\SE_{p_1}\oplus\SE_{p_2}@>q>>\SQ\to 0$$
such that $(\SE_0,\SQ_0)=(\SE,\SQ)|_{\wt X\times\{0\}}=(E,Q)$. The
quotient $\SE_{p_1}\oplus\SE_{p_2}@>q>>\SQ\to 0$ is determined by
the two projections $q_i:\SE_{p_i}\to\SQ$ ($i=1,2$), which can be
choosen to be isomorphisms for $t\neq 0$ since $\SP^{\wt L}$ is
irreduceble. The two maps $q_i$ are given by two matrices
$$\pmatrix
t^{a_1}&0&\hdots&0\\
0&t^{a_2}&\hdots&0\\
\vdots&\vdots&\ddots&\vdots\\
0&0&\hdots&t^{a_r}
\endpmatrix ,\quad \pmatrix
t^{b_1}&0&\hdots&0\\
0&t^{b_2}&\hdots&0\\
\vdots&\vdots&\ddots&\vdots\\
0&0&\hdots&t^{b_r}
\endpmatrix$$
where $0\le a_1\le a_2\le\cdots\le a_r$ and $0\le b_1\le
b_2\le\cdots\le b_r$. When $t=0$, they give the GPB $(E,Q)$. We
recall that when $F$ is not locally free, the numbers $n_1$ and
$n_2$ in the proof of Lemma 1.4 are positive. Thus the two
projections $E_{p_i}\to Q$ are not isomorphism. Namely, there are
$k_1$, $k_2$ such that $a_{k_1}>0$, $b_{k_2}>0$ but $a_j=0$
($j<k_1$) and $b_j=0$ ($j<k_2$). It is clear now that we can
change the positive numbers $a_{k_1}$, ..., $a_r$, $b_{k_2}$, ...,
$b_r$ freely such that the resulted family $(\SE,\SQ)$ has the
property that $(\SE_0,\SQ_0)=(E,Q)$. We modify the $T$-flat
quotient $\SE_{p_1}\oplus\SE_{p_2}@>q>>\SQ\to 0$ by choosing
$a_{k_1}$, ..., $a_r$, $b_{k_2}$, ..., $b_r$ such that
$$\sum^r_{i=k_1} a_i -\sum^r_{i=k_2} b_i=0.$$
Thus we get a $T$-flat sheaf $\SF$ on $X\times T$ such that
$\SF_0=F$. Moreover, on $T\ssm\{0\}$, $\SF$ is obtained from
$\SE|_{\wt X\times(T\ssm\{0\})}$ by identifying $\SE_{p_1}$ and
$\SE_{p_2}$ through the isomorphism
$$q_1\cdot
q_2^{-1}:\SE_{p_1}\to\SE_{p_2}.$$
$\wedge^r\SF|_{X\times(T\ssm\{0\})}$ is obtained from
$det(\SE)|_{\wt X\times(T\ssm\{0\})} =p_{\wt X}^*\wt L$ by
identifying $\wt L_{p_1}\otimes K(T)$ and $\wt L_{p_2}\otimes
K(T)$ through the isomorphism $\wedge^r(q_1\cdot q_2^{-1})$, where
$K(T)$ denote the field of rational functions on $T$. By the
choice of $a_{k_1}$, ..., $a_r$, $b_{k_2}$, ..., $b_r$, we know
that $\wedge^r(q_1\cdot q_2^{-1})$ is the identity map. Thus
$$\wedge^r\SF|_{X\times(T\ssm\{0\})}=(p_X^*L)|_{X\times(T\ssm\{0\})}.$$
\enddemo

\proclaim{Lemma 1.7} For any stable irreducible curve $X$, $\SU_X(r,L)^0$ 
is dense in $\SU_X(r,L)$.\endproclaim

\demo{Proof} Let $\delta$ be the number of nodes of $X$, we will prove
the lemma by induction to $\delta$. When $\delta=1$, it is Lemma 1.6.
Assume that the lemma is true for curves with $\delta-1$ nodes. Then we
show that for any $F\in \SU_X(r,L)$ there is a $T$-flat sheaf $\SF$ on
$X\times T$, where $T=\text{Spec($\Lambda$)}$ and $\Lambda$ is
a discrete valuation ring,
such that\roster
\item $\SF_t=\SF|_{X\times\{t\}}$ is locally free for $t\neq 0$ and
$\SF_0=F$,
\item  $\wedge^r(\SF|_{X\times(T\ssm\{0\})})=p_X^*L$.\endroster

For $F\in \SU_X(r,L)$, we can assume that $F$ is not locally free. Let
$p_0\in X$ be a node at where $F$ is not locally free. Let
$\pi:\wt X\to X$
be the partial normalization at $p_0$ and $\pi^{-1}(p_0)=\{p_1,p_2\}$.
Let $\wt L=\pi^*L$ and $\wt E=\pi^*F/torsion$, then by the same arguments
of Lemma 1.4
$$0\to F@>d>>\pi_*\wt E\to\,_{p_0}\wt Q\to 0.$$
Note that $\wedge^r\wt E=\pi^*(\wedge^rF)/(\text{torsion at $\{p_1,p_2\}$})$
and the cokernel of $\wedge^r\wt E\to\wt L$ at $\{p_1,p_2\}$ is   
$_{p_1}\Bbb C^{n_1}\oplus\,_{p_2}\Bbb C^{n_2}$, we have the morphism
$$\wedge^r\wt E\to \wt L\otimes\SO_{\wt X}(-n_1p_1-n_2p_2)$$
which is an isomorphism outside the nodes of $\wt X$. As the same with proof
of Lemma 1.4, we have the Hecke modification $E$ of $\wt E$ at $p_1$ and $p_2$
such that 
$$0\to \wt E@>h>> E\to\,_{p_1}\wt Q_1\oplus\,_{p_2}\wt Q_2\to 0$$
with $dim(\wt Q_i)=n_i$. Thus $\wedge^rE\cong(\wedge^r\wt E)\otimes\SO_{\wt
X}(n_1p_1+n_2p_2)\to\wt L$ and the generalized parabolic sheaf (GPS) $(E,Q)$ 
defines $F$ by the
exact sequence
$$0\to F@>(\pi_*h)\cdot d>>\pi_*E\to\,_{p_0} Q\to 0,$$
where $Q$ is defined by requiring above sequence exact. The two projections
$E_{p_i}\to Q$ ($i=1,2$) are not isomorphism, thus, by choosing suitable bases of
$E_{p_1}$ and $Q$, they are given by matrices
$$P_1=\pmatrix
1&\hdots&0&\hdots&0\\
\vdots&\ddots&\vdots\\
0&\hdots&1&\hdots&0\\
\vdots&\hdots&\vdots&\ddots\\
0&\hdots&0&\hdots&0
\endpmatrix,\quad P_2=A\cdot\pmatrix
1&\hdots&0&\hdots&0\\
\vdots&\ddots&\vdots\\
0&\hdots&1&\hdots&0\\
\vdots&\hdots&\vdots&\ddots\\
0&\hdots&0&\hdots&0
\endpmatrix\cdot B$$
where $A$, $B$ are invertable $r\times r$ matrices and 
$\text{rank$(P_i)=r_i<r$}$ ($i=1,2$). Since $E\in\SU_{\wt X}(r,\wt L)$, by
the assumption, there is a $T$-flat sheaf $\SE$ on $\wt X\times T$ such that
$\SE_0:=\SE|_{\wt X\times\{0\}}=E$ and $\SE|_{\wt X\times(T\ssm\{0\})}$ locally
free with determinant $p^*_{\wt X}(\wt L)$. Define the morphisms 
$q_i:\SE_{p_i}:=\SE|_{\{p_i\}\times T}\to Q\otimes\SO_T$ ($i=1,2$)
by using matrices
$$Q_1=\pmatrix
1&\hdots&0&\hdots&0\\
\vdots&\ddots&\vdots\\
0&\hdots&1&0\hdots&0\\
0&\hdots&0&t^{a_{r_1+1}}\hdots&0\\
\vdots&\hdots&\vdots&\ddots\\
0&\hdots&0&\hdots \quad c\cdot t^{a_r}
\endpmatrix,\quad Q_2=A\cdot\pmatrix
1&\hdots&0&\hdots&0\\
\vdots&\ddots&\vdots\\
0&\hdots&1&0\hdots&0\\
0&\hdots&0&t^{b_{r_2+1}}\hdots&0\\
\vdots&\hdots&\vdots&\ddots\\
0&\hdots&0&\hdots \quad t^{b_r}
\endpmatrix\cdot B$$
where $t$ is the local parameter of $\Lambda$, $a_{r_1+1}$, ..., $a_r$, $b_{r_2+1}$,
..., $b_r$ are positive integers satisfying 
$a_{r_1+1}+\, \cdots\,+a_r=b_{r_2+1}\,+\cdots\,+b_r$, and $c$ is any constant.
Then these morphisms $q_i$ ($i=1,2$) define a family $(\SE, Q\otimes\SO_T)$ of
GPS, which induces a $T$-flat sheaf $\SF$ on $X\times T$ such that $\SF_0=F$ and
$\SF_t$ ($t\neq 0$) are locally free. The determinant $det(\SF|_{X\times T^0})$,
where $T^0=T\ssm\{0\}$,
is defined by the sheaf $(det(\SE|_{\wt X\times T^0})=p^*_{\wt X}(\wt L)$ 
through the isomorphism
$$det(q_2^{-1}\cdot q_1): (det(\SE|_{\wt X\times T^0})_{p_1}=(\wedge^r\SE_{p_1})|_{T^0}\to
(\wedge^r\SE_{p_2})|_{T^0}=(det(\SE|_{\wt X\times T^0})_{p_2},$$
which is a scale product by $det(Q^{-1}_2\cdot Q_1)=det(AB)^{-1}\cdot c$. Thus we can
choose suitable constant $c$ such that 
$det(\SF|_{X\times T^0})=p^*_X(L)$. We are done.
\enddemo  

\proclaim{Lemma 1.8} When $L$ is not locally free, $\SU_X(r,L)^0$ 
consists of torsion free sheaves of type $1$ at each node of $X$, which
is dense in $\SU_X(r,L)$.\endproclaim

\demo{Proof} The proof follows the same idea. For simiplicity, we assume 
that $X$ has only one node $p_0$. Let $F$ be a torsion free sheaf of rank
$r$ and degree $d$ on $X$ with type $t(F)\ge 1$ at $p_0$. Then
$$deg(\wedge^rF/torsion)=d-t(F)+1.$$
Thus $F\in\SU_X(r,L)^0$ if and only if $t(F)=1$.

For any $F\in\SU_X(r,L)$ of type $t(F)>1$, let $\wt E=\pi^*F/tosion$, then
$$0\to F@>d>>\pi_*\wt E\to\,_{p_0}\wt Q\to 0$$
where $dim(\wt Q)=r-t(F)$. Let $\wt L=\pi^*L/torsion$, then
$deg(\wt L)=d-1$ and $L=\pi_*\wt L$. The condition that $F\in\SU_X(r,L)$
implies that $det(\wt E)=\wt L(-n_1p_1-n_2p_2)$ where $n_i\ge 0$ and
$n_1+n_2=t(F)-1.$ As in the proof Lemma 1.4, let $h:\wt E\to E$ be the Hecke modifications at
$p_1$ and $p_2$ such that $dim(ker(h_{p_1}))=n_1+1$ and
$dim(ker(h_{p_2}))=n_2$. Then we have
$det(E)=det(\wt E)\otimes\SO_{\wt
X}((n_1+1)p_1+n_2p_2)=\wt L(p_1)$, and there is an GPB
$(E,E_{p_1}\oplus E_{p_2}@>q>>Q\to 0)$ such that $\phi(E,Q)=F$, where
$q_i:E_{p_i}\to Q$ ($i=1,2$) satisfy $dim(ker(q_1))=t(F)-n_1-1$ and
$dim(ker(q_2))=t(F)-n_2$.The two projections
$q_i:E_{p_i}\to Q$ ($i=1,2$) are are given by matrices
$$P_1=\pmatrix
1&\hdots&0&\hdots&0\\
\vdots&\ddots&\vdots\\
0&\hdots&1&\hdots&0\\
\vdots&\hdots&\vdots&\ddots\\
0&\hdots&0&\hdots&0
\endpmatrix,\quad P_2=A\cdot\pmatrix
1&\hdots&0&\hdots&0\\
\vdots&\ddots&\vdots\\
0&\hdots&1&\hdots&0\\
\vdots&\hdots&\vdots&\ddots\\
0&\hdots&0&\hdots&0
\endpmatrix\cdot B$$
where $\text{rank$(P_1)=r-t(F)+n_1+1$}$, 
$\text{rank$(P_2)=r-t(F)+n_2$}.$
Let $T=\text{Spec($\Bbb C[t]$)}$ and $\SE=p_{\wt X}^*E$. Choose deformations
$P_i(t)$ of $P_i$ ($i=1,2$) as following
$$\pmatrix
1&\hdots&0&0&\hdots &0&0\\
\vdots&\ddots&\vdots\\
0&\hdots&1&0&\hdots &0&0\\
0&\hdots&0&t&\hdots &0&0\\
\vdots&\ddots&\vdots&\vdots&\ddots&\vdots&\vdots\\
0&\hdots&0&0&\hdots &t&0\\
0&\hdots&0&0&\hdots &0& t
\endpmatrix,\quad A\cdot\pmatrix
1&\hdots&0&0&\hdots&0&0\\
\vdots&\ddots&\vdots\\
0&\hdots&1&0&\hdots&0&0\\
0&\hdots&0&t&\hdots&0&0\\
\vdots&\ddots&\vdots&\vdots&\ddots&\vdots&\vdots\\
0&\hdots&0&0&\hdots&t&0\\
0&\hdots&0&0&\hdots&0&0
\endpmatrix\cdot B$$
where the number of $t$ in $P_2(t)$ is $t(F)-n_2-1$. Then we get
a family $(\SE, Q\otimes\SO_T)$ of GPB on $\wt X\times T$, which
induces a $T$-flat sheaf $\SF$ on $X\times T$ such that $\SF_0=F$
and $\SF_t$ ($t\neq 0$) are torsion free of type $1$. To see that
$\wedge^r\SF_t\cong L$ ($t\neq 0$), we note that 
$det(\SE)=p_{\wt X}^*\wt L(p_1)$ and $L$ is determined by the GPB
$$(det(E)=\wt L(p_1), G\subset det(E)_{p_1}\oplus det(E)_{p_2})$$
where $G$ is the graph of zero map $det(E)_{p_2}\to det(E)_{p_1}$.
Thus we have a non-trivial morphism $\wedge^r\SF_t\to L$, which must be
an isomorphism when $t\neq 0$. 
\enddemo

Next we will prove that $\SU_X(r,L)$ is the underlying scheme of
specialization of the moduli spaces of semistable bundles with
fixed determiant. This in particular implies that 
$\SU_X(r,L)\subset\SU_X(r,d)$ is a
closed subset. Let $S=Spec(R)$ and $R$ be a discrete valuation ring. Let
$\SX\to S$ be a flat proper family of curves with smooth generic
fibre and closed fibre $\SX_0=X$. Let $\SL$ be a relative torsion free
sheaf on
$\SX$ of rank one and (relative) degree $d$ such that $\SL|_X=L$. 
It is well known
that there exists a moduli scheme $f:\SU(r,d)_S\to S$ such that
for any $s\in S$ the fibre $f^{-1}(s)$ is the moduli space of
semistable torsion free sheaves of rank $r$ and degree $d$ on
$\SX_s$ (where $\SX_s$ denote the fibre of $\SX\to S$ at $s$).
Since $\SX$ is smooth over $S^0=S\ssm\{0\}$, there is a family
$\SU(r,\SL|_{S^0})_{S^0}\to S^0$ of moduli spaces of semistable
bundles with fixed determinant $\SL|_{\SX_s}$ on $\SX_s$ ($s\in
S^0$). We have
$$\SU(r,\SL|_{S^0})_{S^0}\subset \SU(r,d)_S.$$
Let $Z$ be the Zariski closure
of $\SU(r,\SL|_{S^0})_{S^0}$ inside $\SU(r,d)_S$. We get a flat family
$$f: Z\to S$$
of projective schemes. For any $0\neq s\in S$, the fibre $Z_s$ is
the moduli space of semistable bundles on $\SX_s$ with fixed
determinant $\SL|_{\SX_s}$.

\proclaim{Lemma 1.9} The fibre $Z_0$ of $f:Z\to S$ at $s=0$ is 
contained in $\SU_X(r,L)$ as a set.\endproclaim

\demo{Proof} We can assume that for any $[F]\in Z_0$
there is a discrete valuation
ring $A$ and $T=Spec(A)\to S$ such that there is a $T$-flat family
of torsion free sheave $\SF$ on $\SX_T=\SX\times_ST\to T$, so that
$$\wedge^r\SF_{\eta}\cong\SL_{\eta},\quad \SF|_{X}\cong F .$$
By Proposition 5.3 of [Se] and its proof (see [Se], it deals with
one node curve, but generalization to our case is straightforward
since its proof is completely local), there is a birational
morphism $\sigma:\Gamma\to \SX_T$ and a vector bundle $\SE$ on
$\Gamma$ such that $\sigma_*\SE=\SF$. Moreover, the morphism $\sigma$
is an isomorphism over $\SX_T\setminus\{x_1,...,x_k\}$.
Since $(\wedge^r\SE)|_{\Gamma_{\eta}}\cong 
(\sigma^*\SL)^{\vee\vee}|_{\Gamma_{\eta}}$, 
note that $(\wedge^r\SE)^{-1}\otimes(\sigma^*\SL)^{\vee\vee}$ is
torsion free (thus $T$-flat), we can extend the isomorphism into
a morphism
$\wedge^r\SE\to (\sigma^*\SL)^{\vee\vee}$. Since $\sigma_*$ and $\sigma^*$
are adjoint functors, $\sigma_*\SO_{\Gamma}=\SO_{\SX_T}$, we have 
$\sigma_*((\sigma^*\SN)^{\vee})=\SN^{\vee}$ for any coherent sheaf $\SN$.
Then, by using 
$\sigma^*(\SN^{\vee})=\sigma^*\sigma_*((\sigma^*\SN)^{\vee})\to
(\sigma^*\SN)^{\vee},$
we have a canonical morphism 
$\sigma_*((\sigma^*\SN)^{\vee\vee})\to \SN^{\vee\vee}.$ In particular, there
is a canonical morphism
$$\sigma_*((\sigma^*\SL)^{\vee\vee})\to \SL^{\vee\vee}\cong\SL$$
which induce
a morphism
$\vartheta:\wedge^r\SF=\wedge^r(\sigma_*\SE)\to\sigma_*\wedge^r\SE\to\SL.$
Modified by some power of the maximal ideal of $A$, we can assume 
the morphism $\vartheta$ being nontrivial on $X$, which means that
$\vartheta$ is an isomorphism on $\SX_T\setminus\{x_1,...,x_k\}$ 
since $X$ is irreducible. Thus $[F]\in \SU_X(r,L).$
\enddemo

\proclaim{Theorem 1.10} $\SU_X(r,L)$ is the closure of $\SU_X(r,L)^0$ in 
$\SU_X(r,d)$. When
smooth curve $\SX_s$ specializes to $\SX_0=X$ and $\SL_s$
specializes to $L$, the moduli spaces $\SU_{\SX_s}(r,\SL_s)$ of
semistable bundles of rank $r$ with fixed determinant $\SL_s$ on
$\SX_s$ specializes to an irreducible scheme $Z_0$ with
$(Z_0)_{\text{red}}\cong\SU_X(r,L)$.\endproclaim

\demo{Proof} Let $\SU(r,d)_S^0\subset\SU(r,d)_S$ be the open subscheme
of torsion free sheaves of type at most $1$. 
Then there is a well-defined $S$-morphism (taking
determinant $det(\bullet)=\wedge^r(\bullet)$)
$$det: \SU(r,d)^0_S \to \SU(1,d)_S.$$
The given family of torsion free sheaves $\SL$ on $\SX$ of rank one
and degree $d$ gives a $S$-point
$[\SL]\in \SU(1,d)_S$. It is clear that
$$Z^0:=(det)^{-1}([\SL])\subset Z$$
and the fibre of $f|_{Z^0}: Z^0\to S$ at $s=0$ is 
irreducible with support $\SU_X(r,L)^0$ (it is also reduced when $L$
is a line bundle). Thus $f^{-1}(0)=Z_0$ contains the closure
$\overline{\SU_X(r,L)^0}$ of $\SU_X(r,L)^0$ in $\SU_X(r,d)$. On the
other hand,
by Lemma 1.9, Lemma 1.8 and Lemma 1.7, we have
$$\overline{\SU_X(r,L)^0}\subset (Z_0)_{\text{red}}\subset\SU_X(r,L)
\subset\overline{\SU_X(r,L)^0}.$$
Hence $\SU_X(r,L)=\overline{\SU_X(r,L)^0}=(Z_0)_{\text{red}}$. In particular,
the fibre of $f:Z\to S$ at $s=0$ is irreducible.
\enddemo

\heading{\S2 Stability and Gieseker functor}\endheading

Let $X$ be a stable curve with $\delta$ nodes $\{x_1,\,...,\,x_{\delta}\}$.
Any semistable curve with stable model $X$ can be obtained from $X$
by destabilizing the nodes $x_i$ with chains $R_i$ ($i=1,...,\delta$)
of projective lines. It will be denoted as $X_{\vec n}$, where
$\vec n=(n_1,...,n_{\delta})$ and $n_i$ is the length of $R_i$ (See [NSe]
for the example of $\delta=1$). Then $X_{\vec n}$ are the curves which
are semi-stably equivalent to $X$,
we use $\pi:X_{\vec n}\to X$ to denote the canonical morphism
contracting $R_1,\,...,\,R_{\delta}$ to $x_1,\,...,\,x_{\delta}$
respectively.
A vector bundle $E$ of rank $r$ on a chain $R=\cup C_i$ of projective
lines is called $positive$ if $a_{ij}\ge 0$ in the decomposition
$E|_{C_i}=\oplus^r_{j=1}\SO(a_{ij})$ for all $i$ and $j$.
A $postive$ $E$ is called $strictly$ $positive$ if for each $C_i$
there is at least one $a_{ij}>0$. $E$ is called $standard$ (resp.
$strictly$ $standard$) if it is positive (resp. strictly positive) and
$a_{ij}\le 1$ for all $i$ and $j$ (See [NSe], [Se]).

For any semistable curve $X_{\vec n}=\cup X^k_{\vec n}$ of genus $g\ge 2$,
let $\omega_{X_{\vec n}}$ be its canonical bundle and
$$\lambda_k=\frac{deg(\omega_{X_{\vec n}}|_{X_{\vec n}^k})}{2g-2},$$
it is easy to see that $\lambda_k=0$ if and only if the irreducible component
$X^k_{\vec n}$ is a component of the chains of projective lines.

\proclaim{Definition 2.1} A sheaf $E$ of constant rank $r$ on $X_{\vec n}$
is called (semi)stable, if for every subsheaf $F\subset E$, we have
$$\chi(F)\,<(\le)\,\frac{\chi(E)}{r}\cdot r(F)\quad
\text{when $r(F)\neq 0,r$},$$
$$\chi(F)\le 0\quad\text{when $r(F)=0$, and $\chi(F)<\chi(E)$ when
$r(F)=r$, $F\neq E$},$$
where, for any sheaf $F$, the rank $r(F)$ is defined to be
$\sum\lambda_k\cdot rank(F|_{X^k_{\vec n}}).$\endproclaim

Let $C=X_{\vec n}$ and
$C_0=X_{(0,n_2,...,n_{\delta})}$ (namly, $C_0$ is obtained from $C$ by
contracting the chain $R_1=\bigcup^{n_1}_{k=1}\Bbb P^1_k$ of projective lines
$\Bbb P^1_k=\Bbb P^1$).

\proclaim{Lemma 2.2} Let $\pi:C\to C_0$ be the canonical morphism,
let $E$ be a torsion free sheaf that is locally free on $R_1$.
If $E|_{R_1}$
is positive and $\pi_*E$ is stable
(semistable) on $C_0$, then $E$ is stable (semistable) on $C$.
In particular, a vector bundle on $X_{\vec n}$
is stable (semistable) if $E|_{R_i}$ ($1\le i\le \delta$) are
positive and
$\pi_*E$ is stable (semistable)
on $X$, where $\pi:X_{\vec n}\to X$ is the canonical morphism
contracting $R_1,\,...,\,R_{\delta}$ to $x_1,\,...,\,x_{\delta}$.
\endproclaim

\demo{Proof} Let $C=\wt C_0\cup R_1$ and $\wt C_0\cap R_1=\{p_1,p_2\}$,
 where $\pi: \wt C_0\to C_0$ is
the partial normalization of $C_0$ at $x_1$. Let $\wt E=E|_{\wt C_0}$,
$E'=E|_{R_1}$. Then we have exact sequence
$$0\to E'(-p_1-p_2)\to E\to\wt E\to 0.\tag2.1$$
If $E|_{R_1}$ is positive and $\pi_*E$ stable (semistable), then
$\pi_*E'(-p_1-p_2)=0.$
For any $E_1\subset E$, consider the sequence (2.1),
let $\wt E_1\subset \wt E$ be the image of $E_1$
in $\wt E$ and $K\subset E'(-p_1-p_2)$ be the kernel of $E_1\to \wt E_1$,
then we have
$$0\to\pi_*E_1\to\pi_*\wt E_1\to R^1\pi_*K=\,_{x_1}H^1(K),$$
and $\chi(E_1)=\chi(\wt E_1)+\chi(K)=\chi(\pi_*\wt E_1)-h^1(K)
\le\chi(\pi_*E_1)$. Since $r(E_1)=r(\pi_*E_1)$,
$$\chi(E_1)-\frac{\chi(E)}{r}r(E_1)\le\chi(\pi_*E_1)-
\frac{\chi(\pi_*E)}{r}r(\pi_*E_1).$$
Thus we will be done if we can check that $\chi(E_1)<\chi(E)$ when
$r(E_1)=r(E)$ and $E_1\neq E$. In this case, the quotient $E_2=E/E_1$
is torsion outside the chains $\{R_i\}$. If $E_2|_R=0$,
where $R=\cup R_i$, then $E_2$ is
a nontrivial torsion and we are done. If $E_2|_R\neq 0$, then
$\chi(E_2)\ge\chi(E_2|_R).$
 Since  $E|_R$
is positive and the surjective map
$$E|_R=\bigoplus^r_{j=1}\SL_j\to E_2|_R\to 0,$$
we have $H^1(E_2|_R)=0$ and there is at least one line bundle $\SL_j$
such that $\SL_j\hookrightarrow E_2|_R$ on a sub-chain. Thus
$\chi(E_2)\ge\chi(E_2|_R)=h^0(E_2|_R)>0$
and $\chi(E_1)<\chi(E)$.
\enddemo

\remark{Remark 2.3} It is easy to show that if $E$ is semistable on
$X_{\vec n}$, then $E$ is $standard$ on the chains and $\pi_*E$ is torsion
free. It is expected that (semi)stability of $E$ also implies the
(semi)stability of $\pi_*E$.\endremark

\proclaim{Definition 2.4} Let $\SC\to S$ be a flat family of stable curves
of genus $g\ge 2$. The associated functor $\SG_S$ (called the
Gieseker functor) is defined as follows:
$$\SG_S: \{ S-schemes\} \to \{sets\},$$
where $\SG_S(T)=\text{set of closed subschemes
$\Delta\subset \SC\times_ST\times_S Gr(m,r)$ such that}$\roster
\item the induced projection map $\Delta\to T\times_S Gr(m,r)$ over $T$
is a
closed embedding over $T$. Let $\SE$ denote the rank $r$ vector bundle
on $\Delta$ which is induced by the tautological rank $r$ quotient bundle
on $Gr(m,r)$.
\item the projection $\Delta\to T$ is a flat family of semistable curves
and the the projection $\Delta\to\SC\times_T T$ over $T$ is the canonical
morphism $\pi:\Delta\to \SC\times_S T$ contracting the chains of projective
lines.
\item the vector bundles $\SE_t=\SE|_{\Delta_t}$ on $\Delta_t$ ($t\in T$) are
of rank $r$ and degree $d=m+r(g-1)$. The qoutients
$\SO^m_{\Delta_t}\to \SE_t$ induce isomorphisms
$$H^0(\SO^m_{\Delta_t})\cong H^0(\SE_t).$$
\endroster
\endproclaim

\proclaim{Lemma 2.5 ([Gi],[NSe],[Se])} The functor $\SG_S$ is represented
by a $PGL(m)$-stable open subscheme $\SY\to S$ of the Hilbert scheme.
The fibres $\SY_s$ ($s\in S$) are reduced, and the singularities of $\SY_s$
are products of normal crossings. A point $y\in\SY_s$ is smooth if and
only if the corresponding curve $\Delta_y$ is a stable curve, namely
all chains in $\Delta_y$ are of length $0$.
\endproclaim

Let $Quot$ be the Quot-scheme of rank $r$ and degree $d$ quotiens
of $\SO^m_{\SC}$ on $\SC\to S$ (we choose the canonical polarization
on any flat family $\SC\to S$ of stable curves of genus $g\ge 2$).
There is a universal quotient
$$\SO^m_{\SC\times_S Quot}\to\SF\to 0$$
on $\SC\times_S Quot\to Quot$. Let $\SR\subset Quot$ be the
$PGL(m)$-stable open subscheme consisting of $q\in Quot$ such that
the quotient map $\SO^m_{\SC\times_S\{q\}}\to \SF_q\to 0$ induces
an isomorphism $H^0(\SO^m_{\SC\times_S\{q\}})\cong H^0(\SF_q)$
(thus $H^1(\SF_q)=0$). We can assume that $d$ is large enough so that
all semistable torsion free sheaves of rank $r$ and degree $d$ on $\SC\to S$
can be realized as points of $\SR$. Let $\SR^s$ ($\SR^{ss}$) be the open
set of stable (semistable) quotients, and let $\SW$ be the closure of
$\SR^{ss}$ in $Quot$. Then there is an ample $PGL(m)$-line bundle
$\SO_{\SW}(1)$ on $\SW$ such that $\SR^s$ (resp. $\SR^{ss}$) is precisely
the set of GIT stable (resp. GIT semistable) points. Thus the moduli
scheme $\SU(r,d)\to S$ is the GIT quotient of $\SR^{ss}\to S$.

Let $\Delta\subset \SC\times_S\SY\times_SGr(m,r)$ be the universal object
of $\SG_S(\SY)$, and
$$\SO_{\Delta}^m\to\SE\to 0$$
be the induced quotient
on $\Delta$ by the universal quotient on Grassmannian over $\SY$. Then
there is a commutative diagram over $S$
$$\CD \Delta @>\pi>>\SC\times_S\SY\\
@VVV     @VVV \\
\SY @=  \SY
\endCD $$

\proclaim{Lemma 2.6} If $S$ is a smooth scheme, then
$\pi_*\SO_{\Delta}=\SO_{\SC\times_S\SY}$
and there is a birational $S$-morphism
$$\theta: \SY \to\SR$$
such that pullback of the universal quotient
$\SO^m_{\SC\times_S \SR}\to\SF\to 0$ (by $id\times\theta$) is
$$\SO_{\SC\times_S\SY}^m\to\pi_*\SE\to 0.$$
\endproclaim
\demo{Proof} Similar with Proposition 6 and Proposition 9 of [NSe]
(See also [Se]).\enddemo

\proclaim{Lemma 2.7} Let $\SY^s=\theta^{-1}(\SR^s)$ and
$\SY^0=\theta^{-1}(\SR^{ss})$. Then
$$\theta: \SY^s\to \SR^s, \quad\theta:\SY^0\to\SR^{ss}$$
are proper birational morphisms.\endproclaim
\demo{Proof} The proof in [NSe] and [Se] for irreducible one node curves
is completely local. Thus can be generalied to general stable curves.
\enddemo

There is a $PGL(m)$-equivariant factorisation (See [NSe], [Se], [Sch])
$$\CD \SY^s@>\imath>>\SY^0@>\imath>> \SH\\
@V\theta VV     @V\theta VV  @V\lambda VV\\
\SR^s@>\imath>>\SR^{ss}@>\imath>>\SW
\endCD$$
and linearisation $\SO_{\SH}(1)$, where $\imath$ is open embedding.
Let $L_a=\lambda^*(\SO_{\SW}(a))\otimes\SO_{\SH}(1)$. Then, for
$a$ large enough, the set $\SH(L_a)^{ss}$ ($\SH(L_a)^s$) of
GIT-semistable (stable) points satisfies:
(i) $\SH(L_a)^{ss}\subset \lambda^{-1}(\SR^{ss})$,
(ii) $\SH(L_a)^s=\lambda^{-1}(\SR^s)$.
By Lemma 2.7, $\theta$ is proper, we have $\lambda^{-1}(\SR^{ss})=\SY^0$
and $\lambda^{-1}(\SR^{ss})=\SY^s$. Thus
$$\SH(L_a)^s=\SY^s=\theta^{-1}(\SR^s),\quad
\SH(L_a)^{ss}\subset\SY^0=\theta^{-1}(\SR^{ss}).$$

\proclaim{Notation 2.8} $\SG(r,d)_S=\SH(L_a)^{ss}//PGL(m)$ is called
(according to [NSe]) the generalized Gieseker semistable moduli space
(or Gieseker space for simplicity). It is intrinsic by
recent work [Sch].\endproclaim

Let $y=(\Delta_y,\SO^m_{\Delta_y}\to\SE_y\to 0)\in\SY^0$. Obviously,
for $y\in \SH(L_a)^{ss}\ssm \SH(L_a)^s$, we have to add extra conditions
besides the semistability of $\pi_*\SE_y$. Alexander Schmitt ([Sch])
recently
figure out a sheaf theoretic condition ($H_3$) (See Definition 2.2.10 in
[Sch])
for $\pi_*\SE_y$, which is
a sufficient and necessary condition for $y\in \SH(L_a)^{ss}$. The
pair $(C,E)$ of a semstable curve $C$ with
a vector bundle $E$
is called $H$-(semi)stable (See [Sch]) if $E$ is $strictly$ $positive$
on the chains of projective lines,
and the direct image (on stable model of $C$) $\pi_*E$
is semistable satisfying the condition ($H_3$).

\proclaim{Theorem 2.9} The projective $S$-scheme $\SG(r,d)_S\to S$
universally corepresents the moduli functor
$\SG(r,d)_S^{\sharp}:\{\text{$S$-schemes}\}\to\{sets\}$,
$$\SG(r,d)_S^{\sharp}(T)=\left\{
\aligned&\text{Equivalence classes of pairs $(\Delta_T,\SE_T)$, where
$\Delta_T\to T$}\\ &\text{is a flat family of semistable curves with
stable model}\\
&\text{$\SC\times _ST\to T$ and $\SE_T$ is an $T$-flat sheaf such that
for}\\
&\text{any
$t\in T$, $(\SE_T)|_{\Delta_t}$ is $H$-(semi)stable vector bundle of}\\&
\text{rank $r$ and degree $d$.}\endaligned\right\}$$
We call that $(\Delta_T,\SE_T)$ is equivalent to $(\Delta'_T,\SE'_T)$
if there is
an $T$-automorphism $g:\Delta_T\to\Delta'_T$, which is identity outside
the chains, such that $\SE_T$ and $g^*\SE'_T$ are fibrewisely
isomorphic.
\endproclaim

\heading{\S3 A Gieseker type degeneration for small rank}
\endheading

Let $\SC\to S$ be a flat family of irreducible stable 
curves and $\SL$ be a line bundle on $\SC$ 
of relative degree $d$.
We simply call the families in $\SG(r,d)_S^{\sharp}(T)$, the
families of semistable Gieseker bundles parametrized by $T$.

\proclaim{Definition 3.1} The subfunctor
$\SG_{\SL}:\{\text{$S$-schemes}\}\to\{sets\}$ of $\SG$ is defined to be
$$\SG_{\SL}(T)=\left\{\aligned&\text{$\Delta\in\SG(T)$ such that for any $t\in T$
there is}\\&\text {a morphism
$det(\SE|_{\Delta_t})\to \pi^*\SL_t$ on $\Delta_t$ which}\\&
\text{is an isomorphism outside the chain of $\Bbb P^1$s}
\endaligned\right\}.$$
\endproclaim

\proclaim{Definition 3.2} The moduli functor
$\SG(r,\SL)_S^{\sharp}$ of semistable
Gieseker bundles with a fixed determinant
is defined to be
$$\SG(r,\SL)_S^{\sharp}(T)=\left\{\aligned&\text{$(\Delta_T,\SE_T)\in
\SG(r,d)_S^{\sharp}(T)$ such that for any $t\in T$}
\\&\text{there exists a morphism
$det(\SE_T|_{\Delta_t})\to \pi^*\SL_t$ on $\Delta_t$}
\\&\text{which is an isomorphism outside the chain of $\Bbb P^1$s}
\endaligned\right\}.$$
When $S=\text{Spec($\Bbb C$)}$, the above defined functor is denoted
by $\SG(r,L)^{\sharp}$.
\endproclaim

Let $S=\text{Spec(D)}$ where $D$ is a discrete valuation ring. Let
$\SC\to S$ be a family of curves with smooth generic fibre and 
closed fibre $\SC_0=X$. Assume that $X$ is irreducible with only one node
$p_0$. Then
we have the following result that is similar with Lemma 1.19 of [Vi].

\proclaim{Lemma 3.3} When $r\le 3$, or $r=4$ but the normalization
$\wt X$ of $X$ is not hyperelliptic, 
the moduli functor
$\SG(r,\SL)_S^{\sharp}$ is a locally closed subfunctor of
$\SG(r,d)_S^{\sharp}$.  More precisely, for any family
$(\Delta_T,\SE_T)\in
\SG(r,d)_S^{\sharp}(T)$, there exists a locally closed subscheme
$T'\subset T$ such that a morphism $T_1\to T$ of schemes factors through
$T_1\to T'\hookrightarrow T$ if and only if
$$(\Delta_T\times_TT_1, pr_1^*\SE_T)\in
\SG(r,\SL)_S^{\sharp}(T_1).$$ Similarly, $\SG_{\SL}$ is a locally closed
subfunctor of $\SG$.\endproclaim

\demo{Proof} Let $\pi:\Delta_T\to \SC\times_S T$ be the birational
morphism contracting the chain of rational curves and $\SL_T$ be 
the pullback $\pi^*\SL$ to
$\Delta_T$. Let $f:\Delta_T\to T$ be the family of semistable
curves (thus $f_*(\SO_{\Delta_T})=\SO_T$). Then the condition that
defines the subfunctor is equivalent to the existence of a
global section of $det(\SE_T|_{\Delta_t})^{-1}\otimes\pi^*\SL_t$
which is nonzero outside the chain $R_t\subset\Delta_t$ of $\Bbb P^1$s.
There is a complex 
$$\SK_T^{\bullet}:\SK_T^0@>\delta_T>>\SK_T^1\tag3.1$$
of locally free sheaves on $T$
such that for any base change
$T_1\to T$ the pullback of $\SK_T^{\bullet}$ to $T_1$
computes the direct image of
$det(\SE_{T_1})^{-1}\otimes\SL_{T_1}$ (which equals to the kernel of
$\delta_{T_1}:\SK^0_{T_1}\to\SK^1_{T_1}$).
There is a canonical closed 
subscheme of $T$ (defined locally by
some minors of $\delta_T$) where $\delta_T$ is not injective.
Replace $T$ by this closed subscheme, we assume that 
$f_*(det(\SE_T)^{-1}\otimes\pi^*\SL)\neq 0.$ Let $U\subset T$ be the
largest open subscheme such that for any $t\in U$ 
$$\text{dim}(H^0(\Delta_t,det(\SE_T|_{\Delta_t})^{-1}\otimes\pi^*\SL_t))=1.$$ 
Let $Y\subset\Delta_U$ be the
support of the cokernel of the map
$$f^*f_*(det(\SE_T)^{-1}\otimes\pi^*\SL)\to det(\SE_T)^{-1}\otimes\pi^*\SL.$$
Let $U_0\subset U$ be the fibre of $U\to S$ at the closed point $0\in S$.
Then $$\pi^{-1}(\{p_0\}\times U_0)\subset \Delta_U$$ 
consists of the chains
of $\Bbb P^1$s. Note that $\pi(Y)\subset X\times U_0$, let
$Y'\subset Y$ be the union of irreducible components $Y_i$ such that
$p_1(\pi(Y_i))\neq p_0$ where $p_1: X\times U_0\to X$ is the projection.
Then we define that $T'=U\setminus f(Y')$.

Let $T_1\to T$ be a morphism. If it factors through 
$T_1\to T'\hookrightarrow T$, it is clear 
$$(\Delta_T\times_TT_1, pr_1^*\SE_T)\in
\SG(r,\SL)_S^{\sharp}(T_1)$$ 
since $t\in T'$ if and only if 
$\text{dim}(H^0(\Delta_t,det(\SE_T|_{\Delta_t})^{-1}\otimes\pi^*\SL_t))=1$
and $$\SO_{\Delta_t}\cong f^*f_*(det(\SE_T)^{-1}\otimes\pi^*\SL)|_{\Delta_t}
\to det(\SE_T|_{\Delta_t})^{-1}\otimes\pi^*\SL_t$$
is an isomorphism outside the chain of $\Bbb P^1$s. On the other hand, if
$$(\Delta_T\times_TT_1, pr_1^*\SE_T)\in
\SG(r,\SL)_S^{\sharp}(T_1),$$
then it factors firstly through the closed subscheme of $T$ where $H^0$
do not vanish. Then we have to show that the image of $T_1$ falls in
the open set $U$, here we need the assumpations that $r\le 3$, or $r=4$
but $\wt X$ is not hyperelliptic. To check it, let $t\in T_1$, then
$\text{dim}(H^0(\Delta_t,det(\SE_T|_{\Delta_t})^{-1}\otimes\pi^*\SL_t))=1$
when $\Delta_t$ has no chain of $\Bbb P^1$s. If $\Delta_t=\wt X\cup R$
has a chain $R$, let $\{p_1,p_2\}=\wt X\cap R$, then
$$H^0(\Delta_t,det(\SE_T|_{\Delta_t})^{-1}\otimes\pi^*\SL_t)
=H^0(\wt X, (det(\SE_T|_{\Delta_t})^{-1}\otimes\pi^*\SL_t)|_{\wt X}
(-p_1-p_2)),$$
which has at most dimension $1$ since
$\text{deg}(det(\SE_T|_{\Delta_t})^{-1}\otimes\pi^*\SL_t)|_{\wt X}(-p_1-p_2))\le 1$
when $r\le 3$, or 
$\text{deg}(det(\SE_T|_{\Delta_t})^{-1}\otimes\pi^*\SL_t)|_{\wt X}(-p_1-p_2))\le 2$
when $r=4$ but $\wt X$ is not hyperelliptic. Thus the morphism $T_1\to T$
factors through $T_1\to U$, then it factors through $T_1\to T'$ by the
definition of functor. 
\enddemo

For simiplicity, we assume that $r$ and $d$ are coprime $(r,d)=1$. In this
case, the functor $\SG(r,d)_S^{\sharp}$ is representable by an irreducible
Cohen-Macaulay $S$-scheme $\SG(r,d)_S\to S$ (See [NSe]), whose fibres are reduced,
irreducible projective schemes with at most normal crossing singularities.
Moreover, there is a canonical proper birational $S$-morphism
$$\theta: \SG(r,d)_S \to \SU(r,d)_S,\tag3.2$$
where $\SU(r,d)_S\to S$ is the family (associated to $\SC\to S$) of moduli
spaces of semistable torsion free sheaves with rank $r$ and degree $d$.

By the above Lemma 3.3, the functor $\SG(r,\SL)_S^{\sharp}$ is representable
by a locally closed subscheme $\SG(r,\SL)_S\subset \SG(r,d)_S$ when $r\le 3$,
or $r=4$ but $\wt X$ is not hyperelliptic.

\proclaim{Lemma 3.4} $\SG(r,\SL)_S\subset \SG(r,d)_S$ is a closed subscheme
of $\SG(r,d)_S$. In fact, for the closed fibre $\SC_0=X$, we have 
$$\SG(r,\SL)_S^{\sharp}(\{0\})=\theta^{-1}(\SU_X(r,\SL_0)).\tag3.3$$
\endproclaim

\demo{Proof} It is enough to prove (3.3). For any 
$(\Delta, E)\in \SG(r,d)_S^{\sharp}(\{0\})$, let 
$$\pi: \Delta\to X$$
be the morphism contracting the chain $R$ of $\Bbb P^1$s. Then, 
by definition of $\theta$,
$$\theta((\Delta,E))=\pi_*(E):=F\in\SU_X(r,d).$$
Note that $F$ has type of $t(F)=\text{deg($E|_R$)}$ (See [NSe]), then
$\pi_*(det(E))$ has torsion of dimension $t(F)-1$ supported at the node
$p_0=\pi(R)$. There is a natural morphism
$$\wedge^rF=\wedge^r(\pi_*E)\to\pi_*(\wedge^rE)=\pi_*(det(E)),$$
which is an isomorphism outside $p_0$. Thus we have an isomorphism
$$\wedge^rF/torsion\cong \pi_*det(E)/torsion$$
since $\text{deg}(\wedge^rF/torsion)=\text{deg}(\pi_*det(E)/torsion)
=d-t(F)+1.$ By using this isomorphism, it is clear that 
$$(\Delta,E)\in \SG(r,\SL)_S^{\sharp}(\{0\})\Longleftrightarrow
\theta((\Delta,E))\in\SU_X(r,d).$$
\enddemo

$\SG(r,\SL)_S$ is in fact a degeneracy loci of a map of vector bundles.
To study it, we recall some standard results (See [FP] for example). 
Let $\varphi: F\to E$ be a morphism of vector bundles on a variety $M$ with $rk(F)=m$ and
$rk(E)=n$. The closed subsets of $M$
$$D_r(\varphi)=\{x\in M\,|\,\text{rank}(\varphi_x)\le r\}$$
are the so called degeneracy locus of $\varphi$. We collect the results into

\proclaim{Lemma 3.5} The codimension of each irreducible component of
$D_r(\varphi)$ is at most
$(n-r)(m-r)$. If $M$ is Cohen-Macaulay and the codimension of each
irreducible of $D_r(\varphi)$
equals to $(n-r)(m-r)$, then $D_r(\varphi)$ is Cohen-Macaulay.
\endproclaim

In (3.1), $rk(\SK^1_T)-rk(\SK^0_T)=g-1$ since $det(\SE_T)\otimes\SL_T$ has
relative degree $0$.
Replace $T$ by an open set $U\subset\SG(r,d)_S$, one sees that
$$\SG(r,\SL)_S=D_{k_0}(\delta_U),\quad k_0=rk(\SK^0_U)-1.$$
In what follows, we will use $\text{Codim($\bullet$)}$ to denote:
codimension of each irreducible component of $\bullet$.
Thus $\text{Codim}(\SG(r,\SL)_S)\le g$, and it is Cohen-Macaulay if
$$\text{Codim}(\SG(r,\SL)_S)=g.$$ 
In particular, let $X$ be the singular
fibre of $\SC\to S$ and $L=\SL|_X$. The closed fibre $G(r,d)$ of
$\SG(r,d)_S\to S$ is the so called generalized Gieseker moduli space
(associated to $X$) of [NSe], which has normal crossing singularities.
The closed fibre of $\SG(r,\SL)_S\to S$, denoted by $G(r,L)$, is the
degeneracy loci $$D_{k_0}(\delta_{U_0})\subset U_0\subset G(r,d)$$
of $\delta_{U_0}:\SK^0_{U_0}\to\SK^1_{U_0}$, where $U_0$ is the closed
fibre of $U\to S$. Thus 
$$\text{Codim}(G(r,L))\le g$$
and $G(r,L)$ is Cohen-Macaulay if $\text{Codim}(G(r,L))=g$.
When $r\le 3$, or $r=4$ but $\wt X$ is not hyperelliptic,
$G(r,L)\subset G(r,d)$ is a closed subscheme that represents 
a moduli functor (See Theorem 3.7 for definition).

\proclaim{Lemma 3.6} $\text{Codim}(G(r,L))=g$. In particular,
$\SG(r,\SL)_S\subset\SG(r,d)_S$ is an irreducible, reduced, 
Cohen-Macaulay subscheme of codimension $g$.
\endproclaim 
\demo{Proof} Assume that $\text{Codim}(G(r,L))=g$. Note that there is
a unique irreducible component of $\SG(r,\SL)_S$ with
codimension $g$ dominates $S$ since
$\SC\to S$ has smooth generic fibre. Thus other irreducible
components (if any) of $\SG(r,\SL)_S$ will fall in $G(r,L)$ and their
codimension in $G(r,d)$ are at most $g-1$ since $\SG(r,d)_S\to S$ is
flat over $S$. This contradicts $\text{Codim}(G(r,L))=g$. Hence 
$\SG(r,\SL)_S\subset\SG(r,d)_S$ is an irreducible, 
Cohen-Macaulay subscheme of codimension $g$. It has to be reduced since
it is Cohen-Macaulay and has a reduced open subscheme. 

Now we prove that $\text{Codim}(G(r,L))=g$ in $G(r,d)$. Let
$J^0_X$ be the Jacobian of line bundles of degree $0$ on $X$. Consider
a morphism $$\phi:G(r,L)\times J^0_X\to G(r,d)$$
that sends any $\{(\Delta,E),\SN\}\in G(r,L)\times J^0_X$  
to $(\Delta,E\otimes\pi^*\SN)\in G(r,d)$, where $\pi:\Delta\to X$ is
the morphism contracting the chain $R$ of $\Bbb P^1$s. We claim
that $$\text{dim}\,\phi^{-1}((\Delta,E_0))\le 1,\quad\text{for any
$(\Delta,E_0)\in G(r,d)$}.$$

Let $\sigma: J^0_X\to J^0_{\wt X}$ be the morphism induced by pulling
back line bundles on $X$ to its normalization $\wt X$. 
The fibres of $\sigma$ are of dimension $1$.
On the other hand, it is easy to see
that the projection $G(r,L)\times J^0_X\to J^0_X$ induces an 
injective morphism
$$\rho:\phi^{-1}((\Delta,E_0))\to J^0_X.$$
To prove the claim, it is enough to show that the image $\text{Im($\rho$)}$
falls in a finite number of fibres of $\sigma$. Note that, for any
$\{(\Delta,E),\SN\}\in\phi^{-1}((\Delta,E_0)),$
we have 
$$det(E)\otimes\pi^*(\SN^{\otimes r})=det(E_0)$$ 
on $\Delta$. Recall that, by definition of $G(r,L)$,
there is a morphism $det(E)\to \pi^*L$ which is an isomorphism outside
the chain $R$ of $\Bbb P^1$s. We have 
$$det(E)|_{\wt X}=\pi^*L|_{\wt X}(-n_1p_1-n_2p_2)=\wt L(-n_1p_1-n_2p_2),$$
where $\wt L$ is the pullback of $L$ to $\wt X$, $n_1$, $n_2$ are nonnegative
integers such that
$$n_1+n_2=\text{deg}(E_0|_R)=t(F_0),\quad F_0:=\pi_*(E_0).$$
Thus
$\sigma\circ\rho(\{(\Delta,E),\SN\})=\sigma(\SN)=\wt\SN\in J^0_{\wt X}$
falls in the set
$$\{\wt\SN\in J^0_{\wt X}\,\,|{\wt\SN}^{\otimes r}=
det(E_0)|_{\wt X}\otimes\wt L^{-1}(n_1p_1+n_2p_2)\},$$
which is clearly a finite set. This proves that fibres of $\phi$ are
at most dimension $1$.

There is a unique irreducible component $G(r.L)^0$ of $G(r,L)$ containing 
$\Delta\cong X$, which has codimension $g$. For any other irreducible
component (if any), say $G(r,L)^+$, all of $\Delta$s in $G(r,L)^+$ must
have chain (with positive length) of $\Bbb P^1$s. Then the image 
$\phi(G(r,L)^+\times J^0_X)$ has to fall in a subvariety of $G(r,d)$, which
has codimension at least $1$. Thus 
$\text{dim$(G(r,L)^+\times J^0_X)\le$}\text{dim$G(r,d)$}$, that is, 
$$\text{Codim$(G(r,L)^+)\ge g$}.$$
By Lemma 3.5, $G(r,L)$ is Cohen-Macaulay of pure codimension $g$.
\enddemo

\proclaim{Theorem 3.7} Let $X$ be an irreducible curve of genus $g\ge 2$
with only one node $p_0$. Let $L$ be a line bundle of degree $d$ on $X$.
Assume that $r\le 3$, or $r=4$ and the normalization of $X$ is not 
hyperelliptic. Then, when $(r,d)=1$, we have
\roster\item There is a Cohen-Macaulay projective scheme $G(r,L)$ of
pure dimension $(r^2-1)(g-1)$, which represents the moduli functor 
$$\SG(r,L)^{\sharp}: \, (\Bbb C-schemes)\to (sets)$$  
which is defined in Definition 3.2.

\item Let $\SC\to S$ be a proper family of curves over a discrete valuation
ring, which has smooth generic fibre $\SC_{\eta}$ 
and closed fibre $\SC_0\cong X$. If there is a line bundle $\SL$ on $\SC$
such that $\SL|_{\SC_0}\cong L$. Then there exists an irreducible, reduced,
Cohen-Macaulay $S$-projective scheme 
$f: G(r,\SL)_S\to S$ such that
$$f^{-1}(0)\cong G(r,L),\quad 
f^{-1}(\eta)\cong\SU_{\SC_{\eta}}(r,\SL_{\eta}).$$
Moreover $G(r,\SL)_S$ represents the
moduli functor $\SG(r,\SL)_S^{\sharp}$ in Definition 3.2.

\item There exists a proper birational $S$-morphism 
$\theta: G(r,\SL)_S\to \SU(r,\SL)_S$ which induces a birational morphism
$\theta: G(r,L)\to \SU_X(r,L).$  
\endroster
\endproclaim


\bigskip





\Refs

\widestnumber\key{KWW}
\ref\key Fa1 \by G. Faltings\paper A proof for the Verlinde
formula \pages 347--374\yr1994\vol 3 \jour J. Algebraic Geometry
\endref

\ref\key Fa2 \by G. Faltings\paper Moduli-stacks for bundles on
semistable curves \pages 489--515\yr1996\vol 304 \jour Math. Ann.
\endref

\ref\key FP\by W. Fulton and P. Pragacz\book
Schubert varieties and degeneracy loci
\bookinfo LNM, 1689
\publaddr Springer-Verlag Berlin-Heidelberg\yr1998\endref

\ref\key Gi \by D. Gieseker\paper A degeneration of the moduli
space of stable bundles \pages 173--206\yr1984\vol 19 \jour J.
Differential Geom.
\endref


\ref\key Ka \by Ivan Kausz\paper A Gieseker type degeneration of moduli
stacks of vector bundles on curves \pages 1--59\jour arXiv:math.AG/0201197
\yr2002
\endref


\ref\key KW \by E. Kunz and R. Waldi\paper Regular differential forms
\vol 79 \jour Contemporary Math.
\endref

\ref\key NR \by M.S.Narasimhan and T.R. Ramadas\paper
Factorisation of generalised theta functions I\pages
565--623\yr1993\vol 114 \jour Invent. Math.\endref

\ref\key NS \by D.S. Nagaraj and C.S. Seshadri\paper Degenerations
of the moduli spaces of vector bundles on curves I\pages 101--137
\yr1997\vol 107 \jour Proc. Indian Acad. Sci.(Math. Sci.)\endref

\ref\key NSe \by D.S. Nagaraj and C.S. Seshadri\paper Degenerations
of the moduli spaces of vector bundles on curves II\pages 165--201
\yr1999\vol 109 \jour Proc. Indian Acad. Sci.(Math. Sci.)\endref

\ref\key Sc \by M. Schessinger\paper Functors of Artin rings
\pages 208--222\yr1968\vol 130 \jour Trans. of AMS. \endref

\ref\key Sch \by A. Schmitt\paper The Hilbert compactification of
the universal moduli space of semistable vector bundles over
smooth curves\jour Preprint\yr2002\endref

\ref\key Se \by C.S. Seshadri\paper Degenerations of the moduli
spaces of vector bundles on curves\vol 1\jour ICTP Lecture Notes
\yr 2000\endref

\ref \key Si \by C. Simpson \paper Moduli of representations of
the fundamental group of a smooth projective variety I \pages
47--129\vol 79 \yr1994\jour I.H.E.S. Publications
Math{\'e}matiques\endref

\ref\key Su1 \by Xiaotao Sun\paper Degeneration of moduli spaces
and generalized theta functions \vol 9\jour J. Algebraic
Geom.\pages 459-527 \yr 2000\endref

\ref\key Su2 \by Xiaotao Sun\paper Degeneration of $SL(n)$-bundles
on a reducible curve\jour 
Proceedings Algebraic Geometry
in East Asia, Japan \yr 2001\endref

\ref\key Te \by M.Teixidor i Bigas\paper Compactifications of
moduli spaces of (semi)stable bundles on singular curves: two
points of view\jour Dedicated to the memory of Fernando Serrano.
Collect Math.\pages 527-548\vol 49\yr1998\endref

\ref\key Vi\by Eckart Viehweg\book
Quasi-projective moduli for polarized manifolds
\bookinfo Ergebnisse der Mathematik und ihrer Grenzgebiete; 3. 
Folge, Bd. 30
\publaddr Springer-Verlag Berlin-Heidelberg\yr1995\endref



\endRefs
\enddocument